%% file: arxivSubmission.tex
\documentclass[12pt]{amsart}
\usepackage{amsfonts,amsmath,amssymb, amsthm,color}
\usepackage{euscript,mathrsfs}

\usepackage[top=1.3in, bottom = 1.3 in, left= 1.1in, right = 1.1in]{geometry}

\usepackage{ulem}

\definecolor{maroon}{RGB}{33,145,140} 
\usepackage[colorlinks=true]{hyperref}
\hypersetup{pdftex, colorlinks=true, linkcolor=maroon, citecolor=maroon, filecolor=blue,urlcolor=blue}
\renewcommand{\em}{\it}

\newcommand{\mbE}{{\mathbb E}}
\newcommand{\HA}{\mathcal{H}}

\newcommand{\mbR}{{\mathbb R}}
\newcommand{\mbW}{{\mathbb W}}
\newcommand{\mbP}{{\mathbb P }}
\newcommand{\mcP}{{\mathcal P }}

\renewcommand{\P}{{\mathcal P}}

\newcommand{\A}{{\mathcal A}}
\newcommand{\D}{{\mathcal D}}
\newcommand{\B}{{\mathcal B}}
\newcommand{\Q}{\mathcal{Q}}

\newcommand{\R}{{\mathcal R}}
\newcommand{\F}{{\mathcal  F}}
\newcommand{\U}{{\mathcal U}}
\newcommand{\mbu}{{\mathbb U}}
\newcommand{\scru}{{\mathscr{U}}}
\newcommand{\ben}{\begin{equation}}
\newcommand{\een}{\end{equation}}
\newcommand{\bena}{\begin{eqnarray}}
\newcommand{\eena}{\end{eqnarray}}
\newcommand{\benas}{\begin{eqnarray*}}
\newcommand{\eenas}{\end{eqnarray*}}

\newtheorem{theorem}{Theorem}[section]
\newtheorem{lemma}[theorem]{Lemma}
\newtheorem{corollary}[theorem]{Corollary}
\newtheorem{proposition}[theorem]{Proposition}
\newtheorem{example}{Example}
\newtheorem{remark}[theorem]{Remark}

\newtheorem{definition}[theorem]{Definition}

\definecolor{vry}{RGB}{253, 231, 37}

\definecolor{vrg}{RGB}{94,201,98}

\definecolor{vrdg}{RGB}{33, 145, 140}

\definecolor{vrb}{RGB}{59,82,139}

\definecolor{vrp}{RGB}{68,1,84}

\definecolor{vro}{RGB}{249,142,9}

\definecolor{vrr}{RGB}{188,55,84}

\definecolor{vrnb}{RGB}{13,8,135}

\def\red#1{\color{\red}{#1}}

\title{Controlled Martingale Problems And Their Markov Mimics }

\author[S. Athreya]{Siva Athreya}
\address{Siva Athreya, International centre for theoretical Sciences, Survey No. 151, Shivakote, Hesaraghatta Hobli, Bengaluru 560089 and Stat-Math Unit, Indian Statistical Institute, 8th Mile, Mysore Road, RVCE P.O., Bengaluru 560059, India}
\email{athreya@icts.res.in}
\thanks{ The work of SA  was supported in
part by Knowledge Exchange grant at ICTS-TIFR}

\author[V.S. Borkar]{Vivek S.\ Borkar}
\address{Vivek S.\ Borkar, Department of Electrical Engineering, Indian Institute of Technology Bombay, Powai, Mumbai 400076, India}
\email{borkar.vs@gmail.com}
\thanks{ The work of VB was supported in part by a S.\  S.\  Bhatnagar Fellowship from the Council of Scientific and Industrial Research, Government of India.
}
\author[N. Gadhiwala]{ Nitya Gadhiwala}
\address{Nitya Gadhiwala, Department of Mathematics, 1984 Mathematics Road, Vancouver, BC Canada V6T 1Z2}
\email{nitya.g20@gmail.com}

\begin{document}
\maketitle 

\begin{abstract}
In this article we prove under suitable assumptions that the marginals of any solution to a relaxed controlled martingale problem on a Polish space $E$ can be mimicked by a Markovian solution of a Markov-relaxed controlled martingale problem. We also show how such `Markov mimics' can be obtained by relative entropy minimisation. We provide many examples where the above results can be applied. 
\end{abstract}

{\em marginals; relative enropy minimisation; stochastic control.} 

\input{main}

\bibliographystyle{plain}
\bibliography{alpha}
\end{document}

%% file: main.tex
\section{Introduction}

Consider a random process $\{X_n\}_{n \geq 0}$, taking values in a Polish space $E$ with the law of $X_0$ given by $\nu$ (say). Consider a Markov  process $\{\tilde{X}_n\}_{n \geq 0}$ on $E$ with the law of $\tilde{X}_0$ given by $\nu$ and transition probabilities specified as $$\mbP(\tilde{X}_{n+1} \in A|\tilde{X}_n = x) := \mbP(X_{n+1} \in A|X_n = x), $$
for all Borel $A \subset E, x \in E, n \geq 0.$  Then it is easy to check using induction that the laws of $X_n$ and $ \tilde{X}_n$ agree for each $n \geq 0$ (In fact, the pair marginals of $(X_n,X_{n+1})$, $(\tilde{X}_n, \tilde{X}_{n+1})$ agree). That is, the Markov process $\tilde{X}$ mimics the one dimensional marginals of the process $X$. We shall say that $\tilde{X}$ is a \textit{Markov mimic} of $X$.  A natural question then is whether this can be done in continuous time. In other words, can one identify conditions under which a  continuous time stochastic processes will have Markov mimics? This issue was taken up by Krylov in \cite{Krylov87} when $X$ is a solution to an Ito differential equation,
$$ X_t =  \int_0^t \beta_s ds +  \int_0^t \delta_s dB_s $$
where $t \geq 0$, $\{B_t\}_{t \geq 0}$ is a standard Brownian motion in ${\mathbb R}^d$ and $\{\delta_t\}_{t \geq 0},\{\beta_t\}_{t\geq 0}$ are bounded processes adapted to its filtration. It was partially addressed by Gy\"{o}ngy \cite{g86} and, in a control theoretic framework, also by Borkar \cite{b86}. They both assume a uniform non-degeneracy condition for the diffusion matrix $\delta_t\delta^T_t$ but  the flavour of the results differ. In \cite{g86}, the existence of a solution to a stochastic differential equation with state-dependent coefficients that mimics the laws of the Ito differential equation is established. The solution, however, need not be Markov unless the stochastic differential equation is well-posed. In \cite{b86}, when  the diffusion matrix is assumed to be a Lipschitz function of state alone, then a  stronger result, viz., that the mimic exists and is a Markov process, is shown.  We shall refer to the process that replicates the one dimensional marginals of a given controlled martingale problem a {\em Markov control mimic} if its controlled extended generator depends on the current state and time alone, without requiring that the process be Markov. We say that it is a {\em Markov mimic} if in addition, the process is Markov. Thus \cite{g86} produces a Markov control mimic whereas \cite{b86} produces a Markov mimic, under their respective sets of assumptions. 

Such results elicited renewed interest following their application in finance, notably due to the work of Dupire \cite{du94}, \cite{du96}. An excellent account of this, along with some important extensions, can be found in \cite{bs13} (see also \cite{f14}). Independently, motivated by stochastic control, there was work by Mikami \cite{m95}, \cite{m99} along similar lines. 

In this paper we address the question of Markov mimics in the very general framework of relaxed controlled martingale problems (see, e.g., Chapter 5 of \cite{ABG} for background and applications). Our  aim is to unify and at the same time extend the existing results. We also point out connections with other results in Markov process theory, controlled or otherwise, by way of remarks. In Theorem \ref{thm:cpmsfm}, under broad assumptions we show the existence of Markov mimics for releaxed controlled martingale problems and point out its implications in stochastic control. We also show that  discounted occupational measures can be mimicked by time homogeneous Markov processes (see Theorem \ref{thm:main2}).

Our assumptions guarantee existence of Markov controls and a Markov solution (see Remark \ref{rem:ass}). It is trivial to note that if there is no Markov solution to the martingale problem then the problem of finding a Markov mimic is vacuous. Theorem \ref{thm:cpmsfm} shows existence of Markov mimics and also  has implications in stochastic controls where costs (that are to be optimised) depend on one-dimensional marginals (see Remark \ref{rem:thmscm}).  
Examples where our broad assumptions hold are discussed in Section \ref{sec:examples}. 

More recently, a renewed interest in this topic was generated by optimal transport, wherein minimization of  entropy (to be precise, \textit{relative entropy}, i.e.\  Kullback-Leibler divergence) as a route to Markov mimics was explored in \cite{l13}, \cite{ac20}, \cite{bF20}. Stochastic control problems are closely related to Schr\"odinger bridges and the Monge-Kantorovich optimal transport problems. See \cite{cgtp21} for a survey for understanding connections between one-time marginal flows in  control problems with McCann displacement in optimal transport. See also \cite{l13} for a survey of the Schr\"odinger problem and its connections to optimal transport. In \cite{ac20}, the authors consider a generalisation of the Schr\"odinger problem, namely the so called Br\"odinger problem. The objective is to minimise relative entropy, with respect to a base measure $\mbP_0$, over a set of measures $\mbP$ with certain prescribed constraints on the marginals.  Under markovian assumption on the base measure $\mbP_0$, it is shown that if the optimisation problem has a unique solution then it is also Markov \cite[Theorem 4.1]{ac20}. 

In Section \ref{sec:relent}, we focus on lowering relative entropy. In Proposition \ref{prop:emimic} we show that if $\mbP$ is non-Markov then the relative entropy can be lowered  by a suitable markovianisation procedure of $\mbP$ (see Definition \ref{def:mtion})  that preserves marginals. The proof is adapted from \cite{b91} and a similar technique is used in \cite{ac20} as well. We also give a sufficient condition on the constraint set for existence of such Markov mimics that minimize relative entropy, in Theorem \ref{thm:scm} and Corollary \ref{cor1}.
See Remark \ref{rem:ent} for a discussion on the uniform integrability assumption imposed  in the hypothesis of the two results and also how Corollary \ref{cor1} may  be used for an alternative method of Markov selection. Some related literature is as follows. 
 In \cite{cgtp22}, the authors consider a formulation of minimisation of relative entropy for diffusion with killing. The unbalanced optimal transport problem is handled via suitable augmentation, see \cite[Problem 7]{cgtp22}, where our results are also applicable. See also \cite{m21} for results on stochastic control with fixed marginals and connections to Schr\"odinger bridges and optimal transport. 

The rest of the article is organised as follows. In the next section we introduce the controlled martingale problem, required assumptions and prove our  main results. By adapting the argument of \cite{b86}, we show existence of  Markov mimics (Theorem \ref{thm:cpmsfm}) and that discounted occupational measures can be mimicked by time homogeneous Markov processes (see Theorem \ref{thm:main2}). In Section \ref{sec:examples}, we give representative examples that illustrate the applicability of our main result.  Section \ref{sec:relent} develops the alternative approach of entropy minimisation in significant generality, see Proposition \ref{prop:emimic} and Theorem \ref{thm:scm}. We conclude the paper with  a couple of examples that illustrate the applications of our results on entropy minimisation to questions in optimal transport.


\section{Markov mimics}

Let $\mbu$ be a Polish space. For a generic Polish space  $E$, $\P(E)$ will denote the space of probability measures on $E$ endowed with the Prokhorov topology and ${\mathcal B_b(E)}$ will denote the space of bounded measurable functions $E \rightarrow \mathbb{R}$. Let $D([0,\infty), ;E)$ be the Polish space of r.c.l.l.\  paths from $[0,\infty) \rightarrow E$ with the Skorokhod topology and let $\scru$ be the space of all measurable maps from $[0,\infty) \rightarrow \P(\mbu)$.

\noindent {\em Topology on $\scru$:} Let $\{f_i\}$ be a countable dense set in the unit ball of $C(\overline{\mbu})$, where $\overline{\mbu}$ is the standard compactification of $\mbu$, i.e., the closure of its usual homeomorphic embedding into $[0,1]^\infty$. Then $\{f_i\}$ is a convergence determining class for $\P(\mbu)$. For $U \in \scru$, let
$$ \alpha_i(t) := \int_\mbu f_i(u)U_t(du), \qquad i = 1,2, \ldots$$
Then $\alpha_i$ has measurable paths and $|\alpha_i(t)| \leq 1$ for all $t \geq 0$. For $T >0$, let ${\mathcal B}_T$ denote the space of measurable maps $[0,T] \rightarrow [-1,1]$ with the weak$^{\star}$-topology of $L^2[0,T]$ relativized to it. Let ${\mathcal B}$ denote the space of measurable maps $[0,\infty) \rightarrow [-1,1]$  with the corresponding inductive topology, i.e., the coarsest topology that renders continuous the map ${\mathcal B} \rightarrow {\mathcal B}_T$ that maps $x \in {\mathcal B}$ to its restriction to ${\mathcal B}_T$, for every $T >0$. Let ${\mathcal B}^\infty$ be the countable product of ${\mathcal B}$ with the product topology.

Next, note that the map $\phi :\P(\mbu) \rightarrow [-1,1]^\infty$ defined by
$$ \mu \in \scru \mapsto \left ( \int f_1 d\mu,  \int f_2 d\mu, \int f_3 d\mu, \ldots\right) \in \mathcal{B}^\infty$$
is continuous, one-to-one with a compact domain, and hence is a homeomorphism onto its range. Equivalently, we denote the map $ \P(\mbu) \rightarrow (\alpha_1,\alpha_2,\ldots) \in {\mathcal B}^\infty$ for $\alpha_i(\cdot) := \int f_id\mu(\cdot)$ also as $\phi : \scru \rightarrow {\mathcal B}^\infty$.  We relativize the topology of ${\mathcal B}$ to $\phi(\scru)$
and topologize ${\scru}$ with the coarsest topology that renders 
$\phi:\scru \rightarrow \phi(\scru)$ a homeomorphism.

From \cite[Theorem 2.3.2]{ABG} it is immediate that $\scru$ is compact and metrizable, hence Polish. Furthermore, \cite[Theorem 2.3.3]{ABG} also implies that if $U^n \rightarrow U$ in $\scru$ as $ n \rightarrow \infty$ and $f \in C([0,T] \times \mbu)$ for some $T >0,$ then
$$\int_0^T \int_\mbu f(t,u) U^n_t(du)dt \longrightarrow \int_0^T \int_\mbu f(t,u) U_t(du)dt$$ as $ n \rightarrow \infty.$

\begin{definition} Let $\A$ be a linear operator with domain $\D(\A) \subset {\mathcal B}_b(E)$ and range $\R(\A) \subset {\mathcal B}_b(E)$. Let $\nu \in \P(E).$ An $ E$ valued process $\{X_t : t \geq 0\}$   on a probability space $(\Omega, \F, \mbP)$ is said to be a solution to the  martingale problem for $(\A, \nu)$ with respect to a filtration $\{\F_t\}_{t \geq 0}$ if \\

\begin{enumerate}
  \item[(i)] $X$ is progressively measurable with respect to $\{\F_t\}_{t \geq 0}$ ,

  \item[(ii)] the law of $X_0$ is $\nu$, and,

  \item[(iii)] for all $f \in \D(\A)$ 
  \begin{equation} \label{eq:mp}
  f(X_t) - f(X_0)- \int_0^t \A f(X_s)ds
  \end{equation}
  is an $\{\F_t\}_{t \geq 0}$ martingale under $\mbP$.
  \end{enumerate}
\end{definition}

\begin{definition}\label{def:CMPstrong}
Let $\A$ be a linear operator with domain $\D(\A) \subset {\mathcal B}_b(E)$ and range $\R(\A) \subset {\mathcal B}_b(E \times \mbu)$. Let $\nu \in \P(E).$ An $ E \times \mbu$ valued process $\{(X_t,U_t) : t \geq 0\}$   on a probability space $(\Omega, \F, \mbP)$ is said to be a solution to the controlled martingale problem for $(\A, \nu)$ with respect to a filtration $\{\F_t\}_{t \geq 0}$ if : 

\begin{enumerate}
  \item[(i)] $(X,U)$ is progressively measurable with respect to $\{\F_t\}_{t \geq 0}$,

  \item[(ii)] the law of $X_0$ is $\nu$, and,

  \item[(iii)] for all $f \in \D(\A)$ 
  \begin{equation} \label{eq:srmp}
  f(X_t) - f(X_0)- \int_0^t \A f(X_s, U_s)ds
  \end{equation}
  is an $\{\F_t\}_{t \geq 0}$ martingale under $\mbP$.
  \end{enumerate}
  Correspondingly an $E \times \P(\mbu)$-valued process $(X,U)$ defined on a probability space $(\Omega, {\mathcal F}, \mbP)$ is said to be a solution to the {\em relaxed controlled martingale problem} for $(\A, \nu)$ with respect to a filtration $\{\F_t\}_{t \geq 0}$ if 
  \begin{enumerate}
  \item[(i)] $(X,U)$ is progressively measurable with respect to $\{\F_t\}_{t \geq 0}$, 

  \item[(ii)] the law of $X_0$ is $\nu$, and,

  \item[(iii)] for all $f \in \D(\A)$ 
  \begin{equation} \label{eq:rmp}
  f(X_t) - f(X_0)- \int_0^t \int_{\U}\A f(X_s, u) U_s(du)ds
  \end{equation}
  is an $\{\F_t\}_{t \geq 0}$ martingale under $\mbP$.
  \end{enumerate}
\end{definition} 

We define the operator $\bar{\A} : {\mathcal D}({\mathcal A}) \rightarrow {\mathcal B}_b(E \times \P(\mbu))$ by
$$\bar{\A}f(x ,\nu) = \int_\mbu {\A}f(x,u)\nu(du), $$
where $ f \in {\mathcal D}({\mathcal A}), x \in E$ and $\nu \in \P(\mbu).$ Consequently, \eqref{eq:rmp}  can be written as 
$$f(X_t) -f(X_0)-\int_0^t \bar{\A}f(X_s, U_s)ds. $$

For $u \in \mbu$ and $\nu \in \P(\mbu)$, we shall define ${\A}^u: {\mathcal D}({\mathcal A}) \rightarrow C_b(E \times \mbu)$ and $\bar{\A}^\nu: {\mathcal D}({\mathcal A}) \rightarrow C_b(\P(\mbu \times E))$ respectively, by
$$\A^uf(x) = \A f(x,u) \mbox{ and } \bar{\A}^\nu f(x)= \bar{\A} f(x,\nu),$$
where $ f \in {\mathcal D}({\mathcal A}), x \in E, u \in \mbu$ and $\nu \in \P(\mbu).$

Finally, we say that a  (relaxed) controlled martingale problem  has a {\em Markov solution} if it has a solution $(X_t, U_t)$ for each initial condition $x \in E$ such that the collection $x \mapsto \mbP_x :=$ the law of this solution for $x \in E$, satisfies the Chapman-Kolmogorov equation. \\

 We make the following assumptions.

 {\bf Assumptions:}

 \begin{enumerate}
\item[(A1)] There exists a countable set $\{g_k\} \subset {\mathcal D}(\A)$ such that 
$$ \{ (g, {\mathcal A}g): g \in {\mathcal D}(\A) \subset \mbox{bp-closure}  \{(g_k, {\A}g_k): k \geq 1\},$$
where bp-closure is bounded pointwise closure (see \cite[ Section 3.4, page 111]{EK} for a definition).

\item[(A2)]  ${\D}(\A)$  is an algebra that separates points in $E$ and contains constant functions. Also, $\A {\mathbf 1} = 0$ where ${\mathbf 1}$ is the constant function identically equal to $1$.
\item[(A3)] Given a $u \in \mbu$ and $x \in E$,
\begin{enumerate}
\item[(i)] there exists a r.c.l.l.\,solution to the martingale problem for $(\A^u, \delta_x)$, with $\delta_x$ being the Dirac measure at $x \in E$,
\item[(ii)]$\A^u$ is dissipative (see \cite[Section 1.2, page 11]{EK} for a definition),
\item[(iii)]${\mathcal D}(\A^u)$ is dense in ${\mathcal B}_b(E),$ 
\item[(iv)] $\R(I-A) = {\mathcal B}_b(E)$.
\end{enumerate}
 \item[(A4)] By a standard measurable selection theorem (see \cite[Lemma 1]{benes70}), given a solution $(X_\cdot, U_\cdot)$ to the  relaxed controlled martingale problem for $({\mathcal A}, \nu)$,  there exists a measurable map $\bar{v}:   E\times [0,T] \rightarrow \P(\mbu)$ such that
\begin{equation} \label{eq:rmc}
{ {\mbE[\bar{\A}f(X_s, U_s)| X_s] = \bar{\A}^{\bar{v}(X_s,s)}f(X_s) }}\end{equation}
a.s. for $s \in [0,T].$  Assume that the martingale problem for $(\bar{\A}^{\bar{v}(X_s, s)}, \nu)$  with $\{\mathcal{F}_t\}$ replaced by $\{\mathcal{F}^X_t\} :=$ the natural filtration of $X$, has a solution for all $\nu$. 

\item[(A5)]Suppose we are given a solution $(X_\cdot, U_\cdot)$ to the  relaxed controlled martingale problem for $({\mathcal A}, \nu)$ and $\bar{v}(\cdot,\cdot)$  is defined as in (A4). Then the relaxed control problem for $(\bar{\A}^{\bar{v}(X_s,s)},\delta_x)$, $x \in E$,  has a Markov solution.


\end{enumerate}

Before we proceed, we make a few observations concerning the above assumptions.
\begin{remark} $\mbox{}$ \label{rem:ass}
{\rm
\begin{itemize}


  \item {\em Markov and Stationary Markov Controls: } A control of the form $U_t = v(X_t,t)$, with $t\geq 0$ and measurable $v: E\times\mbR_+ \to \mbu$ is said to be a Markov control. It is said to be a stationary Markov control if it is of the form $U_t = v(X_t),$ for $t\geq 0$ and a  measurable $v: E \to \mbu$.  Analogous definitions apply for relaxed controls. If $X$ is Markov, resp.\ time-homogeneous Markov, the control may be taken to be a Markov, resp.\ stationary Markov control. This can be proved along the lines of  \cite[Theorem 2.2.23, p.\ 46]{ABG}.  The converse, however, is not true, as, e.g., in case of uncontrolled degenerate diffusions with bounded continuous coefficients. In fact, the entire set of solutions can be characterized in this case as  in  \cite[Section 12.3]{SV}. \\


  \item {\em Existence of a Markov solution: } 
If the martingale problem for $(\bar{\A}^{\bar{v}(X_s,s)}, \nu)$ is well-posed for all $\nu$, then the additional requirement of the existence of a Markov solution is automatically satisfied. In general, if the set of solution measures for each initial condition $x$ is nonempty and compact in $\mathcal{P}(D([0,\infty);E))$, then a procedure due to Krylov \cite{Krylov} (see  also \cite[Section 12.2]{SV}) yields a Markov selection for diffusion.  For an alternative selection procedure for finite dimensional uncontrolled diffusions, see \cite{BS10,BS23}.  See \cite[Theorem 5.9]{EK} for sufficient criteria when  Krylov's Markov selection can be done on solutions of martingale problems.
In our case this is true, e.g., if the relaxed control $U_\cdot$ is of the form 
\begin{equation} \label{eq:msol}
U_t = F(X([0,t]),t) 
\end{equation}
 for a continuous $F(\cdot,t): D([0,t], E) \rightarrow \mathcal{P}(\mbu)$ for each $t \geq 0$. This would imply that $\bar{v}(\cdot, t)$ is continuous for each $t$ and facilitate the desired compactness of solution measures (usually proved by first establishing tightness and therefore relative compactness  thereof using standard criteria and then showing that each subsequential limit is a legitimate solution measure, for which continuity of coefficients plays a role). 
 In general one can only guarantee that $\bar{v}(\cdot,t)$ are measurable. For finite dimensional uncontrolled diffusions, in this generality, one can possibly use a stochastic differential inclusion \cite{Kisiel} that will facilitate a compact set of admissible laws for a given initial condition. We flag this as a direction for future research.\\

  \item {\em Weak vs strong solutions: } 
   In general, the martingale problem for $(\A^{v(X_s,s)}, \nu)$ (alternatively, $(\bar{\A}^{\bar{v}(X_s,s)}, \nu)$ ) has to be interpreted in the weak sense, i.e., the underlying probability space is not specified a priori, but only the existence thereof is asserted, and uniqueness is interpreted in terms of uniqueness of the laws.  We shall refer to this as the `weak formulation' to distinguish it from the `strong formulation' of Definition \ref{def:CMPstrong}. The following, however, holds:{}\\

 If $(X,U)$ is a solution to the martingale problem $(A^u, \nu)$ with respect to $\{\mathcal{F}_t\}$, then $(X,U')$ with $U'_t = v(X_t,t) \ \forall t \geq 0$, is a solution to the martingale problem $(\A^{v(X_s,s)}, \nu)$ with respect to $\{\mathcal{F}_t\}$. Conversely, if the latter problem has a (weak) solution on some probability space, then one can construct a copy in law of $(X,U)$ on a possibly augmented version of this probability space.  This follows as in  \cite[Theorem 2.3.4, p.\ 52]{ABG}, where this result is proved for controlled diffusions.\\

\end{itemize}
}
 \end{remark}

We are now ready to state our main result.\\

\begin{theorem} \label{thm:cpmsfm}
Assume (A1)-(A5). Given any solution $(X, U)$ to a relaxed controlled martingale problem for $({\mathcal A}, \nu)$, there exists a Markov control $\bar{v}(\cdot, \cdot)$ and a solution $X^\prime$ to the relaxed controlled martingale problem for $({\mathcal A}, \nu)$ with this Markov control and $X_0 \stackrel{d}{=} X^\prime_0$, such that $X, X^\prime$ have identical one dimensional marginals. Furthermore, $X^\prime$ can be taken to be a Markov solution.
\end{theorem}

\begin{proof}
By assumption (A1)-(A3), let $\{(X_t,U_t)\}_{t \geq 0}$ be a solution to the controlled martingale problem.  In view of (A4) and (A5), let $\tilde{X}$ be a Markov solution to the martingale problem for $(\bar{\A}^{\bar{v}}, \delta_x)$, where $\bar{\nu}$ is as in \eqref{eq:rmc}. 

Fix $t>0$.  Using (A1)-(A3), we will choose a version of $\bar{v}(X_\cdot,\cdot)$ such that the (two parameter) transition semigroup $T_{s,t}$ with $0 \leq s \leq t$ is   a strong contraction semigroup with generator $\bar{\A}^{\bar{v}}$ on the Banach space of bounded measurable functions with supremum norm (see \cite[Phillips - Lumer Theorem, p.\ 250]{yos95}). Then 
by {{\cite[Proposition 1.1.5]{EK} we have for $f \in \mathcal{D}(\A)$, $T_{s,t}f \in \mathcal{D}(\A)$ and }}
\begin{equation}\label{eq:gen}
\frac{\partial}{\partial s}T_{s,t}f +{\bar{\A}}^{\bar{v}(\cdot,s)}T_{s,t}f = 0, 
\end{equation}
for all $s \in [0,t]$. As $(X_t,U_t)$ satisfies the controlled martingale problem,  applying \eqref{eq:mp} to the function $T_{s,t}f$ we have that
$$T_{t,t}f(X_t) - T_{0,t}f(X_0)-\int_0^t \left(\frac{\partial}{\partial s}T_{s,t}f(X_s)+\bar{\A} T_{s,t}f(X_s, U_s)\right)ds
$$
is an $\{\F_t\}_{t \geq s}$ martingale under $\mbP$.  From the above, we then have
\begin{eqnarray*}
 \mbE[f(X_t)] - \mbE[f(\tilde{X}_t)] &=&\mbE[T_{t,t}f(X_t)] - \mbE[T_{0,t}f({X}_0)] \\
 &=&\mbE \left[\int_0^t \left[\frac{\partial}{\partial s}T_{s,t}f(X_s)+\bar{\A} T_{s,t}f(X_s, U_s)\right]ds \right].\\
\end{eqnarray*}
Using \eqref{eq:gen} leads to 
\begin{align*}
\mbE[f(X_t)] -\mbE[f(\tilde{X}_t)] 
&= \int_0^t  \mbE[-{{{\bar{\A}}^{\bar{v}(\tilde{X}_s,s)}}}T_{s,t}f(X_s)+ \bar{\A} T_{s,t}f(X_s,U_s)]ds\\
&= \int_0^t  \mbE[ \mbE[-{{{\bar{\A}}^{\bar{v}(\tilde{X}_s,s)}}}T_{s,t}f(X_s)+ \bar{\A} T_{s,t}f(X_s,U_s) \mid X_s]]ds\\
&=\int_0^t  \mbE[ -{{{\bar{\A}}^{\bar{v}(\tilde{X}_s,s)}}}T_{s,t}f(X_s) + \mbE[ \bar{\A} T_{s,t}f(X_s,U_s) \mid X_s]]ds\\
&=0
\end{align*}
by \eqref{eq:rmc}. {{As $ f \in  \mathcal{D}({\bar{\A}})$ was arbitrary, by (A1), (A2) we have that $X$ and $\tilde{X}$ have the same marginals.}}
\end{proof}

\begin{remark} $\mbox{}$ \label{rem:thmscm}
\begin{itemize}

\item A weaker version of this result appears as  \cite[Theorem 2.4, p.\ 1552]{BB}, where it is proved that the one dimensional marginals can be mimicked by a process controlled by a Markov control. But it is not asserted or claimed that the latter process itself is Markov. A similar observation applies to \cite[Theorem 4.6, p.\ 516]{g86} where under non-degeneracy condition on the diffusion matrix  of an Ito differential equation, the one-dimensional marginals are mimicked by a stochastic differential equation with measurable drift and diffusion matrix. The latter conditions ensure only existence and not uniqueness (\cite[Section 2.6]{Krylov08}). Consequently it can have non-Markov solutions. This issue is avoided in \cite{b86} by means of a stronger additional condition on the diffusion matrix, viz., that it is a Lipschitz function of the current value of the process alone. Then the Markov controlled process is itself Markov.\\

\item One immediate implication for stochastic control problems wherein the cost or reward depends only on one dimensional marginals, is that the existence of an optimal non-anticipative control implies the existence of an optimal Markov control.\\

\item Theorem \ref{thm:cpmsfm} says that, given a solution $(X_\cdot, U_\cdot)$ to the relaxed controlled martingale problem, there exists a measurable map $\bar{v}: E \times [0,T] \rightarrow \P(\mbu)$ such that
\begin{equation} \label{eq:mc}
{{\mbE[\A f(X_s, U_s)| X_s] = \A^{\bar{v}(X_s,s)}f(X_s) }}\end{equation}
a.s.\ for $s \in [0,T],$ and the    martingale problem for $(\A^{\bar{v}(X_s,s)}, \nu)$ has a solution $\widetilde{X}(\cdot)$ that is a Markov process with the same one dimensional marginals as $X_{\cdot}$. If all solutions of the latter martingale problem have identical marginals for \textit{every} choice of the initial distribution $\nu$, then by Theorem 4.4.2, p.\ 184, of \cite{EK}, the martingale problem in fact has a unique solution. 
\end{itemize}

 \end{remark}
 
 We conclude this section with a related result that the so called \textit{$\alpha$-discounted occupation measure} $\mu$  for $(X_{\cdot},U_{\cdot})$ defined by
$$\int f(x,u)\mu(dx,du) := \alpha\mbE\left[\int_0^\infty e^{-\alpha t}f(X_t,U_t)dt\right], \ f \in C_b(E\times\A),$$
can be replicated by a Markov mimic controlled by a \textit{stationary} Markov relaxed control $\bar{v}(\cdot)$. \\

\begin{theorem}\label{thm:main2} Assume (A1)-(A5). Given a discount factor $\alpha > 0$ and any solution $(X,U)$ of the relaxed controlled martingale problem for $({\mathcal A}, \nu)$, there exists a relaxed stationary Markov control $\bar{v}(\cdot)$ and a solution $X'$ to the relaxed controlled martingale problem for $({\mathcal A}, \nu)$ with this relaxed stationary Markov control with $X_0\stackrel{d}{=}X_0'$, such that $X,X'$ have identical marginals and therefore identical $\alpha$-discounted occupation measures.

\end{theorem}

\begin{proof}

Let $E$ be a Polish space,  $\mbu$ be a compact metric space.  Let $(X, U)$ be a relaxed controlled martingale problem for $(\A, \nu)$ satisfying (A1)-(A5).

 Let $x \in E$ and $\alpha >0$. Define a probability measure $\mu$ on $E \times \mbu$ by 
 \begin{equation}
  \int f d\mu  = \alpha \mbE\left [ \int_0^\infty e^{-\alpha s} f(X_s, U_s)ds \right],
 \end{equation}
for any bounded continuous function $f: E \times \mbu \rightarrow \mbR$. Let the marginal on $E$ of $\mu$ be denoted by $\eta$ and let $\bar{v}(du\mid y)$ denote the conditional distribution of $u$ given $y$ under $\mu.$ Let $X^\prime$ be the solution to  the relaxed Markov controlled martingale
problem with $(\bar{\A}^{\bar{v}},\nu)$. 

Let $k \in {\mathcal D}(\bar{\A})$ and $\psi(x) = \mbE_x \left [ \int_0^\infty e^{-\alpha s} k(X^\prime_s, \bar{v}(X^\prime_s))ds \right]$ $:=\mbE \left [ \int_0^\infty e^{-\alpha s} k(X^\prime_s, \bar{v}(X^\prime_s))ds \mid X_0 =x\right] $. Then from the definition of $\bar{v}(\cdot)$, it follows that
$$ \bar{\A}^{\bar{\nu}} \psi(x) - \alpha \psi(x) +k(x,\bar{v}(x))=0$$
and 
\begin{align*}
\mbE \left [ \int_0^\infty e^{-\alpha s} k(X_s, \bar{v}(X_s))ds \right]
&=\mbE\left [ \int_0^\infty e^{-\alpha s} \left(- \bar{\A}^{\bar{\nu}} \psi(X_s) + \alpha \psi (X_s) \right)ds \right]\\
&= \lim_{T\uparrow\infty}\mbE \left [ \int_0^T e^{-\alpha s} \left(- \bar{\A}^{\bar{\nu}} \psi(X_s) + \alpha \psi (X_s) \right)ds \right]\\
&= \int \psi(x) \nu(dx) - \lim_{T\uparrow\infty}e^{-\alpha T}\mbE\left[\psi(X_T)\right] = \int \psi(x) \nu(dx).
\end{align*}
This establishes the claim.
\end{proof}

\begin{remark} \label{rem:thm2}
As in \cite{b86}, one can consider stationary Markov controls and try to mimic laws at exit times.  Suppose $(X, U)$ is a solution  to a relaxed controlled martingale problem on $\mbR^d$. Suppose $X^\prime$ is a time-homogeneous Markov solution to the relaxed controlled martingale problem with a stationary Markov control, $\bar{r}(X^\prime_t) \in \P(\mbu)$ say. Let $\tau = \inf \{ t \geq 0 : X^\prime_t \not \in D\},$ with $D \subset \mbR^d$ and $\mbE[\tau] < \infty.$  Then one could imitate the arguments in Theorem \ref{thm:main2} and show that  if $X_0 \stackrel{d}{=} X^\prime_0$, then $X_\tau \stackrel{d}{=} X^\prime_\tau$. 
Define for $f \in C^2_0(\R^d)$, $h : \mbR^d \rightarrow \R$  by $$h(x) = \mbE [ f(X^\prime_{\tau^\prime}) \lvert X^\prime_0 = x)]. $$

Then we will require that $h$ solves the Dirichlet problem given by
 \begin{equation}\label{eq:elliptic}
 {\bar{\A}}^{\bar{r}(x)} h = 0, \qquad h = f \mbox{ on } \partial D \end{equation}
Then as in \cite{b86}, an application of Dynkin's formula  yields the result. One would need to impose additional assumptions on $\bar{\A}$ so that \eqref{eq:elliptic} holds. 

\end{remark}

\subsection{Examples} \label{sec:examples}

In this section we shall discuss several examples where Theorem \ref{thm:cpmsfm} is applicable.  We discuss controlled martingale problems that arise naturally in applications, satisfying the hypotheses of Theorem \ref{thm:cpmsfm} and Theorem \ref{thm:main2}. First we note that, if the problem is well-posed, i.e.\ the respective martingale problem has a unique solution, then the solution is already Markov.  

We begin with an example from finite dimensional diffusions.

\begin{example} \label{eg:diffusion} 
Let $E= \mbR^d$, $\mathbb{U}$ be any compact metric space, $\nu \in \P(E)$ and $S_d$ denote the set of all symmetric non-negative definite $d\times d$ real matrices.
For $i,j\in\{1,2,\ldots,d\}$, define  $a_{ij}:E\rightarrow \mbR$ and $b_i:E\times \mathbb{U} \rightarrow\mbR$ such that $a_{ij}$ and $b_i$ are bounded and measurable for all $i,j$ and $a = [a_{ij}]\in S_d$. Let $\A$  be a linear operator with $\D(\A) = C^2_0(\mbR^d)$ be given by
  \begin{equation} \label{eq:cdiff}
  \A f(x,u) = \sum_{i=1}^d b_i(x,u)\frac{\partial f}{\partial x_i}(x) + \frac{1}{2} \sum_{i,j=1}^d a_{ij}(x)\frac{\partial^2 f}{\partial x_ix_j}(x).  
  \end{equation}
As $\D(\A) = C^2_0(\mathbb{R}^d)$, it is easy to see that (A1) and (A2), (A3) (ii), (iii), (iv) are satisfied. Then, by \cite[Theorem 6.1.7]{SV} there is a solution  $(X,U)$ to the martingale problem associated to $(\A^u,\delta_x)$, so (A3) (i) holds.  By \cite[Theorem 4.1]{kurtz98} or \cite[Theorem 2.4]{BB}, (A4) holds. Finally, from \cite[Theorem 12.2.3]{SV} or \cite[Theorem 5.19]{EK},  (A5) holds  when $U$ is as in \eqref{eq:msol} and $b,a$ are bounded continuous functions. 
\end{example}

Next we consider the case of pure jump diffusion.

\begin{example} \label{eg:jump}
Let $E$ be a locally compact Polish space. Let $\mbu$ be a Polish space. Let $\lambda : E \times \mbu \rightarrow \mbR_+$ be a non-negative, measurable functions bounded on compact sets. Let  $\gamma(x,u,A)$ be a transition function on $E \times \mbu \times {\mathcal B}(E).$ Let $(X_t,U_t)$ be a solution to the controlled martingale problem for $(\B,\nu)$ where
$$ \B f(x,u) =  \lambda(x,u)\int_{\mbR^d} [f(x+y) - f(x) ]\gamma(x,u,dy)
$$
for $f \in C_0(E^\Delta)$, where $E^\Delta$ is a one point compactification of $E.$ Further assume that for $x\in E, u \in \mbu$ and for $f \in C_0(E)$,  $$  \lambda(x,u)\int_{\mbR^d} \mid f(x+y) - f(x) \mid \gamma(x,u,dy) < \infty. $$
As $\D(\B) = C_0(E)$, it is easy to see that (A1), (A2), and  (A3) (ii), (iii) and (iv) are satisfied. From \cite[Exercise 15 in p.\ 263]{EK} or  \cite[Example 3.5]{kurtz98} there is a solution  $(X,U)$ to the martingale problem associated to $(\A^u,\delta_x)$, so (A3) (i) holds.  By \cite[Theorem 4.1]{kurtz98} or \cite[Theorem 2.4]{BB}, (A4) holds. Finally, from \cite[Theorem 5.19]{EK},  (A5) holds  when $U$ is as in \eqref{eq:msol} and $\lambda$ is a bounded continuous function. 
\end{example}

By \cite[Example 3.3]{kurtz98}, for $E = \mathbb{R}^d$, a  linear combination of $\A$ as in Example \ref{eg:diffusion} and $\B$ as in Example \ref{eg:jump} will also satisfy (A1)-(A4). (A5) will also follow if they both satisfy the respective hypothesis required in each of the examples.  We now present an example in the infinite dimensional setting.

\begin{example}\label{eg:infinitedim}
Let $E=H$ be a real separable Hilbert space. Let $U$ be a closed unit ball $K$ of another real separable Hilbert space, with the weak topology. Let  $F:H \rightarrow K$ be continuous, $B:K \rightarrow H$ be bounded linear,  $W(\cdot)$ be an $H$-valued Wiener process with covariance given by a trace class operator $Q$.  Let $L$ be an  infinitesimal generator of a differentiable compact semigroup of contractions $\{S(t)\}_{t \geq 0}$ on $H$ such that $L^{-1}$ is a bounded self-adjoint operator with discrete spectrum. Let $\{e_i : i \geq 1\}$ be a
 CONS in  $H$ such that they are  eigenfunctions of $L^{-1}$ with corresponding eigenvalues $\{\lambda_i^{-1}: i \geq 1\}$. Let $P_n: H \rightarrow \mbR^n$ be the map defined by $P_n(x) = [\langle x,e_1 \rangle, \ldots, \langle x,e_n \rangle].$
Let ${\mathcal D}(\A) = \{f \circ P_n:  f\in C_0^2(\mbR^n), n \geq 1\} \subset C_b(H).$ Define $A: {\mathcal D}(\A) \rightarrow C_b(H \times U)$ by 
$$ \A(f\circ P_n)(h,u)= \sum_{i=1}^n \langle e_i, (F(h) + Bu-\lambda_i h) \rangle \frac{\partial f}{\partial x_i} \circ P_n(h)+ \frac{1}{2} \sum_{i,j=1}^n \langle e_i, Q e_j \rangle  \frac{\partial^2 f}{\partial x_i \partial x_j}.$$
By definition of $\D(\A)$, it is easy to see that (A1) and (A2), (A3) (ii) and (iii) are satisfied. From 
\cite[Theorem 8.1]{pz14}  (A3) (i) holds.  \cite[Example 3 and Theorem 2.4]{BB},  ensure that (A4) holds.
 \end{example}
We conclude this section with an example from nonlinear filtering theory. This arises from control problems for diffusion with partial observations.
\begin{example}  \label{eg:filtering}
$E=\mcP(\mbR^d)$ and $\mbu$ be any compact metric space. Let $\A$ be as in Example \ref{eg:diffusion}. Let $${\mathcal D}(B) = \left \{ f \in C_b(\mcP(\mbR^d)) : \begin{array}{l} f(\mu) = g(\int f_1d\mu, \int f_2 d\mu, \ldots, \int f_n d\mu), \mu \in \mcP(\mbR^d), \\ \\ \mbox{ for some } n \geq 1, g \in C_0^2(\mbR^n), f_1,  \ldots, f_n \in {\mathcal D}(A) \end{array} \right\}.$$ Let $B$ be a linear operator from $B: {\mathcal D}(B) \rightarrow C_b(\mcP(\mbR^d) \times \mbu)$ defined by
\begin{align*}
&Af(\mu,u) = \sum_{i=1}^n \frac{\partial g}{\partial x_i} 
(\int f_1d\mu, \ldots, \int f_n d\mu) \int Af_i(\cdot, u) d\mu\\
&+ \frac{1}{2} \sum_{i,j=1}^n  \frac{\partial^2 g}{\partial x_i \partial x_j}(\int f_1d\mu, \ldots, \int f_n d\mu) \times \\ & \hspace{1in}\left \langle \int f_ih d\mu-\int f_i d \mu \int h d\mu, \int f_jh d\mu-\int f_j d\mu \int h d\mu \right\rangle
\end{align*}
By definition of $\D(B)$ and the hypotheses assumed on $\A$ from Example \ref{eg:diffusion}, it is easy to see that (A1) and (A2), (A3) (ii), (iii) and (iv) are satisfied. From discussion in 
\cite[Section 8.2, 8.3]{ABG}  there is a solution  $(X,U$ to the martingale problem associated to $(\A^u,\delta_x)$, so (A3) (i) holds.  Finally, \cite[Example 4, Theorem 2.4]{BB} ensure  that  (A4) holds. Such a treatment is also possible for stochastic evolution equations (see \cite{lz21}).
\end{example}

\section{Minimizing Relative Entropy} \label{sec:relent}
 We begin with the definition of relative entropy between two probability measures on  Polish spaces.
 \begin{definition}
 For a  Polish space $S$ endowed with its Borel $\sigma$-field $\B(S)$, let $\mbP, \mbP_0$ be probability measures on $(E, \B(S))$ with $\mbP << \mbP_0$. Let $\mbE[ \ \cdot \ ], E_0[ \ \cdot \ ]$ denote the respective expectation operators and let $\Lambda := \frac{d\mbP}{d\mbP_0}$ denote the Radon-Nikodym derivative of $\mbP$ w.r.t.\ $\mbP_0$. We define the relative entropy (equivalently, the Kullback-Leibler divergence) of $\mbP$ with respect to $\mbP_0$ as
$$D(\mbP\|\mbP_0) := \mbE_0[\Lambda\log\Lambda] = E[\log\Lambda].$$
     
 \end{definition}
Let $T>0$, $E$ be a polish space and  let $D([0,T];E)$  be the polish space of r.c.l.l.\ paths in $E$ with the Skorokhod topology. Let $\mbP_0$ denote a reference probability measure on $D([0,T];E)$ under which the coordinate process is Markov. Let  $\{X_t\}_{0 \leq t \leq T} \in D([0,T];E)$ be a r.c.l.l.\ process  whose law  $\mbP$ satisfies: $\mbP_t :=$ the restriction of $\mbP$ to  $D([0,t];E)$ is absolutely continuous w.r.t.\ $\mbP_{0t} :=$ the restriction of $\mbP_0$ to $D([0,t]; E)$, $\forall \ t \geq 0$. For $t > 0$, let $\Lambda_t := \frac{d\mbP_{t}}{d\mbP_{0t}}$  with $\Lambda_0 := 1$, and for $0 < s < t$, let
$$\tilde{\Lambda}_{s,t} := \frac{\Lambda_t}{\Lambda_s}$$
 when $\Lambda_s > 0$ and $= 0$ otherwise.
\begin{definition} Suppose for $0 \leq u < v $ if ${\mathcal{F}}_{u,v} = \sigma({X}_r : u \leq r \leq v)$. Let $0 < s < t \leq T$.
We say that $s$ is a  `\textit{Markov point}' for $X$ if $\mathcal{F}_{s,t}$ and $\mathcal{F}_{0,s}$ are conditionally independent given $X_s$. \end{definition}

{\em Markovianisation:} Fix $0<s <t \leq T$. For $T\geq 0$, given a process $X \in D([0,T];E)$, we shall use the notation $X([s,t])$ for $0 \leq s < t \leq T$  to denote the restriction of $X$ to $[s,t]$, viewed as an element of $D([s,t];E)$.  Construct $\breve{X}$ on the path space  $D([0,t];E)$  as follows. The process $\breve{X}[0,s]$ has  the same law as that of $X([0,s])$. Let the conditional law of $\breve{X}([s,t])$, given $\breve{X}([0, s)$ be the  conditional law of $X([s,t])$ given $X_s$. Then the values of $\breve{X}([0,s])$ and $\breve{X}([s,t])$ are matched at $s$ and the concatenation thereof can be viewed as an element of $D([0,t];E)$.  
 More precisely, for any $n \geq1$, $1 \leq i \leq {n-1}$, $0 \leq s_i \leq s_{i+1} \leq s \leq t_i \leq t_{i+1}  \leq t$ and Borel $A_i, B_i, B$ 
\begin{align*} 
&\mbP(\breve{X}_{s_i} \in A_i, 1 \leq i \leq n) : = \mbP({X}_{s_i} \in A_i, 1 \leq i \leq n),\\ 
& \mbP(\breve{X}_{t_i} \in B_i, 1 \leq i \leq n) : = \mbP({X}_{t_i} \in B_i, 1 \leq i \leq n),\\
&\mbP(\breve{X}_{s} \in B) : = \mbP({X}_{s} \in B), \\
&\mbP(\breve{X}_{t_i} \in A_i, 1 \leq i \leq n \mid  \breve{X}_s \in B):=
 \mbP({X}_{t_i} \in A_i, 1 \leq i \leq n \mid  {X}_s \in B) \mbox{ and } \\ 
&\mbP(\breve{X}_{t_i} \in A_i, 1 \leq i \leq n \mid \breve{X}_{s_i} \in B_i,1 \leq i \leq n \mbox{ and } \breve{X}_s \in B)\\
&\hspace{1in}:=
 \mbP({X}_{t_i} \in A_i, 1 \leq i \leq n \mid  {X}_s \in B).
\end{align*}
 Note that the process $\breve{X}$ is well defined on the canonical path space $D([0,t]; E)$ by the above definition.  In particular, the process  $\breve{X}$ has law which is identical to that of $X$ on the interval $[0,s]$ and $[s,t].$ Suppose for $0 \leq u < v $ if $\breve{\mathcal{F}}_{u,v} = \sigma(\breve{X}_r : u \leq r \leq v)$. Then for any $s \leq r \leq t$ the conditional law of $\breve{X}_r$ given $\breve{\mathcal F}_{0,s}$ is the same as the conditional law of $\breve{X}_r$ given $\breve{X}_s$.  This ensures that  $\breve{\mathcal F}_{0,s}$ and  $\breve{\mathcal F}_{s,t}$ are conditionally independent given $\breve{X}_s.$  Note that  the law of 
$\breve{X}$, $X$ will not be the same in $[0,t]$ unless $s$ is a \textit{Markov point} for $X$.
\begin{definition} \label{def:mtion}
The process  $\breve{X}$ constructed above will be defined  as the `\textit{markovianizing}' of $X$ at time $s$. A set $\mathcal{G}$ of probability measures on $D([0,T];E)$ is said to be closed under markovianization at (say) $s\in [0,T]$ if the law of $\breve{X}$ above is in $\mathcal{G}$ whenever the law of $X$ is.
\end{definition}

\begin{lemma} \label{lem:conent} Let  $\{X_t\}_{0 \leq t \leq T}$, $\mbP$, $\mbP_t $, $\mbP_0$, $\mbP_{0t} $ be as above. Fix $0<s <t\leq T$ and  $\{\breve{X}_u\}_{ 0 \leq  u \leq T}$ be the markovianisation of $X$ at $s$ and $\breve{P}$ denote its law.  Let $\breve{P}_t$  be the restriction of $\breve{P}$ to $D([0,t]; E)$, $\forall \ 0\leq t \leq T$. Then $\breve{P}_t \ll \mbP_{0t}$ and Radon-Nikodym derivative $\breve{\Lambda}_t$ of $\breve{P}_{t}$ w.r.t.\ $\mbP_{0t}$ is given by $$\mbE_0\left [ {\Lambda_t} \bigg | { \mathcal{F}_{s,t}}\right].$$ 
\end{lemma}
\begin{proof}
{Let $g\in C_b(D([0,s];E)), h \in C_b(D([s,t];E))$. Let $\breve{E}[ \ \ ]$ denote the expectation under $\breve{P}$. Then by construction of $\breve{X}$ we have 
\begin{eqnarray} \label{eq:mas}
\breve{E}[g(\breve{X}([0,s]))h(\breve{X}([s,t]))] &=&\breve{E}[g(\breve{X}([0,s]))\breve{E}[h(\breve{X}([s,t]))|\breve{\mathcal{F}}_{0,s}]] \nonumber \\
&=& \breve{E}[g(\breve{X}([0,s]))\breve{E}[h(\breve{X}([s,t]))|\breve{X}_s]]
\end{eqnarray}
Now as  the conditional law of $X([s,t])$ given $X_s$ is the same as the  conditional law of $\breve{X}([s,t])$ given $\breve{X}(s)$ and the  law of $X([0,s])$ is the same as the law of  $\breve{X}([0,s])$ we have  that  
\begin{eqnarray} \label{eq:XsBX}
\breve{E}[g(\breve{X}([0,s]))\breve{E}[h(\breve{X}([s,t]))|\breve{X}(s)]] &=& \mbE[g(X([0,s]))\mbE[h(X([s,t]))|X_s]]. 
\end{eqnarray}
Recall that $\Lambda_t$ is the Radon-Nikodym derivative of $\mbP_t$ w.r.t $\mbP_{0,t}$.  So,
\begin{eqnarray} \label{eq:ens1}
\mbE[g(X([0,s]))\mbE[h( X([s,t]))|X_s]] &=&
\mbE_0\left[\Lambda_s g(X([0,s])) \mbE[h(X([s,t]))|X_s \right]\nonumber\\ &=&\mbE_0\left[g(X([0,s]))\Lambda_s {\left(\frac{\mbE_0[h(X([s,t]))\Lambda_t|X_s]}{\mbE_0[\Lambda_t|X_s]}\right)}\right]. \nonumber\\
\end{eqnarray}
Here the second equality follows from the change of measure formula for conditional expectations. 
As  $\sigma(X_s) \subset \mathcal{F}_{s,t}$ and $X([s,t])$ is $\mathcal{F}_{s,t}$ measurable we have 
\begin{equation}
E_0[h(X([s,t])\Lambda_t|X_s]  = E_0[E_0[h(X([s,t])\Lambda_t|\mathcal{F}_{s,t}]|X_s] \\
= E_0[h(X([s,t])E_0[\Lambda_t|\mathcal{F}_{s,t}]|X_s] 
\end{equation}
Therefore, using this in \eqref{eq:ens1} we have
\begin{eqnarray}\label{eq:ens1a}
\mbE[g(X([0,s]))\mbE[h( X([s,t]))|X_s]] &=&\mbE_0\left[g(X([0,s]))\Lambda_s\left({\frac{\mbE_0[h(X([s,t]))\mbE_0[\Lambda_t|\mathcal{F}_{s,t}]|X_s]}{\mbE_0[\Lambda_t|X_s]} }\right)\right] \nonumber \\
\end{eqnarray}
We know that the coordinate process is Markov under $\mbP_0$, so the above equals  
 \begin{eqnarray} \label{eq:ens2}
&& \mbE_0\left[g(X([0,s]))\Lambda_s\left(\frac{\mbE_0[h(X([s,t]))\mbE_0[\Lambda_t|\mathcal{F}_{s,t}]|\mathcal{F}_{0,s})]}{\mbE_0[\Lambda_t|\mathcal{F}_{0,s}]}\right)\right]  \nonumber\\
&=&\mbE_0\left[g(X([0,s]))\Lambda_s\left(\frac{\mbE_0[h(X([s,t]))\mbE_0[\Lambda_t|\mathcal{F}_{s,t}]|\mathcal{F}_{0,s}]}{\Lambda_s}\right)\right] \nonumber\\
&=&\mbE_0[g(X([0,s]))\mbE_0[h(X([s,t]))\mbE_0[\Lambda_t|\mathcal{F}_{s,t})]|\mathcal{F}_{0,s})]]
\end{eqnarray}
where the second last line follows from  the martingale property of $\Lambda_t$ under $\mbP_0$. Thus from \eqref{eq:mas}, \eqref{eq:XsBX}, \eqref{eq:ens1a}, and \eqref{eq:ens2} we have
\begin{equation} \label{eq:rnpr}
\breve{E}[g(\breve{X}([0,s]))h(\breve{X}([s,t]))]= \mbE_0\bigg [g(X([0,s]))h(X([s,t])) \,\,\mbE_0[\Lambda_t|X([s,t])]\bigg]. 
\end{equation}
It is easy to see that \eqref{eq:rnpr} holds for functions of the form $(x,y) \mapsto \sum_{i=1}^n \alpha_i g_i (x) h_i(y)$ with
$\alpha_i \in \mbR, x \in D([0,s];E)), y \in C_b(D([s,t]), g_i\in C_b(D([0,s];E))$ and\\ $ h_i \in C_b(D([s,t];E)).$
The claim follows via an application of Stone-Weierstrass theorem.
}
\end{proof}

\begin{proposition} \label{prop:emimic}
  {{Let  $\{X_t\}_{t \geq 0}$, $\mbP$, $\mbP_T $, $\mbP_0$, $\mbP_{0T} $ be as above. If there exists an $s\in(0,T]$ that is not a Markov point for $X$, then the process $\{\breve{X}_t\}_{t \geq0}$ obtained by markovianising $X$ at $s$ satisfies
\begin{enumerate}
	\item the marginals of $\breve{X}$ is the same as the marginals of $X$ and,
	\item  the law of $\breve{X}$ given by the probability measure $\breve{P}$ on $D([0,\infty);E)$  with  $\breve{P}_T :=$ the restriction of $\breve{P}$ to  $\D([0,T];E)$, satisfies 
	\begin{equation}  D(\breve{P}_T || \mbP_{0T}) <  D(\mbP_T || \mbP_{0T}), \ \end{equation}
	for all $T >s.$ 
\end{enumerate}
  }}
\end{proposition} 

\begin{proof} 
{Let $\breve{\Lambda}$ be the Radon-Nikodym of $\breve{P}_t$ w.r.t. $\mbP_{0t}$. Using Lemma \ref{lem:conent}, we have }
\begin{eqnarray*}
\mbE_0[\breve{\Lambda}_t\log\breve{\Lambda}_t] &=&  \mbE_0[\mbE_0[\Lambda_t| \mathcal{F}_{s,t}]\log(\mbE_0[\Lambda_t| \mathcal{F}_{s,t})]\\
&<&  {E_0[\Lambda_t \log(\Lambda_t)],} 
\end{eqnarray*}
 where the last line follows by the conditional Jensen's inequality and the strong convexity of the map $x \in [0,\infty) \mapsto x\log x \in \mathbb{R}. $
{The proposition readily follows from this, the construction of $\breve{X},$ and the definition of relative entropy.}
\end{proof}

From Proposition \ref{prop:emimic}, it follows that among all $E$-valued r.c.l.l.\ processes that have  the same marginals as $X$ and have laws absolutely continuous with respect to $\mbP_0$, the minimiser of relative entropy, if one exists, is Markov. A more general claim holds :

\begin{theorem} \label{thm:scm}
Let $\mbP_0$ be a reference probability measure on $D([0,T];E)$ under which the coordinate process is Markov.
Suppose  $\Q \subset \P(D([0,T];E))$ is a set of probability measures absolutely continuous w.r.t.\ $\mbP_0$  that is { closed under markovianization }at any $s \in [0,T]$, and let $\HA := \{\Lambda_T  := \frac{d\mbP_T}{d\mbP_{0T}} | \ \mbP \in  \Q\}$ equipped with  $\sigma(L_1,L_\infty)$ topology ($:=$ the weak topology on $L_1$). 
\begin{enumerate}
\item[(a)]If $\mbE_0[\Lambda_T\log\Lambda_T]$ attains its minimum on $\HA$, then  the minimizer is unique and is the law of a Markov process. 

\item[(b)]Suppose that  there exists $h: \mbR_+ \rightarrow \mbR_+$ such that

\begin{equation} \label{eq:h}
\lim\limits_{ t \rightarrow \infty} \frac{h(t)}{t \log(t)} = \infty
\end{equation} 
and 
\begin{equation}
\sup_\HA \mbE_0\left[h(\Lambda_T)\right] < \infty, \label{eq:dagger}
\end{equation} 
 then $\mbE_0[\Lambda_T\log\Lambda_T]$ attains its minimum on $\HA$.
\end{enumerate}
\end{theorem}

\begin{proof} 
Define $f: \HA \mapsto [0,\infty)$ by $$f(\Lambda_T) = \Lambda_T\log \Lambda_T.$$
The map $x \mapsto x\log x$ is convex and continuous on $\mbR_+$. Hence 
$$f(\Lambda_T) = \sup_{g\in\mathcal{C}}g(\Lambda_T),$$
where 
\begin{eqnarray*}
\mathcal{C} &:=& \{g: \mbR_+ \to \mbR : g(x) = ax + b \ \ \mbox{for some } a, b \in \mbR \ \ \mbox{and} \ \ g(x) \leq f(x) \ \forall \ x \geq 0\}. \nonumber
\end{eqnarray*}
Therefore $f(\cdot)$ is lower semi-continuous a.s.\ on $\HA$ and  hence so is the function $F: \HA \mapsto \mbR_+$ given by $F(\Lambda_T) := \mbE_0[\Lambda_T\log\Lambda_T]$. Consequently, if $F$ attains its minimum on $\HA$,  there is a unique minimizer due to the {strong convexity of $f$.} Suppose that the minimiser is not a Markov process. {Then it has a non-Markov point $s$ and the $\breve{\Lambda}$ (:= the corresponding Radon-Nikodym derivative  for the probability  measure $\breve{P}_t$ as defined in the proof of Proposition \ref{prop:emimic}, w.r.t $\mbP_{0t}$) will have a strictly lower value of $F$, a contradiction.} Hence the unique minimiser is a Markov process.

 Under (\ref{eq:dagger}), $\HA$ is uniformly integrable by the de la Vall\'{e}e Poussin theorem (\cite{DMNH}, p.\  24II). Therefore it is relatively compact and relatively sequentially compact in the $\sigma(L_1,L_\infty)$ topology by the Dunford-Pettis compactness criterion  (\cite{DMNH}, p.\  27II). It is also easy to check that $\HA$ is closed. Therefore $f$ attains its minimum on $\HA$ by the Weierstrass theorem.
\end{proof}

\begin{corollary}\label{cor1}
Let $\Q_m :=$ the closed subset of $\Q$ whose elements have the same one dimensional marginals as some prescribed element of $\Q$ at $t \in$ some $\mathcal{T} \subset [0,\infty)$. Define $\HA_m$ correspondingly, analogously to the $\HA$ above. Then a unique minimiser  of $F$ on $\HA_m$ exists and will be a Markov process.
\end{corollary}

This is immediate on observing that $\Q_m$ is closed under markovianization at any point.\\

\begin{remark} $\mbox{}$ \label{rem:ent}
\begin{itemize}
\item Note that \eqref{eq:h} is satisfied  for e.g. $h(x) = x^{1+\epsilon} $ for any $\epsilon >0.$ Furthermore, under \eqref{eq:dagger}, Theorem \ref{thm:cpmsfm} then shows the existence of a minimizer.  This result could be of independent interest.
  
\item A priori, the law of a Markov process whose marginals match those of a given random process need not be absolutely continuous with respect to the law of the latter. For example let $\{B_t\}_{t \in [0,1]}$ be a Brownian motion  and define the process to be
$$Y_t = \left \{ \begin{array}{cc}
     B_t& t \in [0, \frac{1}{2}), \\
     B^\prime_t& t \in [\frac{1}{2},1],
     \end{array} \right. $$ where $\{B^\prime_t\}_{\frac{1}{2} \leq t \leq 1}$ is an independent Brownian motion such that $B^\prime_\frac{1}{2}$ is Normal with mean $0$ and variance $\frac{1}{2}.$  The above result then shows that under (\ref{eq:dagger}), there is at least one Markov mimic for which absolute continuity holds.

\item There is one case where the uniform integrability of $\HA_m$ is easy to obtain without a condition such as (\ref{eq:dagger}). Note that $\Lambda_t, t \geq 0$, is a multiplicative functional of the sample path, which makes $\log\Lambda_t, t \geq 0$, an additive  functional. There are cases (e.g., diffusion processes) where $E_0[\Lambda_t\log\Lambda_t] = E[\log\Lambda_t]$ in fact depends only on one dimensional marginals of the process. In this case, this quantity is a constant on $\HA_m$.  Uniform integrability is often easy to check in these scenarios. Even in some situations where this additive functional does not depend only on one-dimensional marginals, this may give an easy route for verifying uniform integrability, e.g., for reflected diffusions where the additive functional involves local time at the boundary. 

\item One interesting result about controlled martingale problems of the type studied in Theorem \ref{thm:cpmsfm} is as follows. Fix an initial distribution $\nu$. Define an equivalence relation, denoted by `$\approx$' between two solution processes for the controlled martingale problem for a prescribed controlled extended generator $\A^u$ as follows: Set $(X_\cdot, u_\cdot) \approx (X'_\cdot,u'_\cdot)$ if their one dimensional marginals agree Lebesgue-a.e. 
The following is proved in \cite{ABG}, see Theorem 6.4.16, p.\ 241, \cite{ABG}, extending an earlier result for controlled diffusions from \cite{b93}. 

\begin{theorem} The extreme points of the closed convex set (in quotient topology)  of such equivalence classes are singletons containing a Markov process.\end{theorem}

 Corollary \ref{cor1} now gives us, under the additional hypotheses of absolute continuity w.r.t.\ a common reference measure, an additional piece of information, viz.\  that every equivalence class contains a Markov solution as well.
 \end{itemize}
\end{remark}

To  illustrate the application of Proposition \ref{prop:emimic} and Theorem \ref{thm:scm} we begin with an example of finite dimensional diffusion discussed in Section \ref{sec:examples}.

\begin{example}[Example \ref{eg:diffusion} contd.] In \cite{lzks21}, trajectories of cellular development are modelled using 
\begin{equation}  \label{eq:wsde}
  dX_t = -\nabla \Psi (t,X_t) +\sigma dB_t,
\end{equation}
with $X_t$ taking values in a compact smooth Riemannian manifold without boundary $E$, $\Psi: [0,1] \times E \rightarrow \mbR$ is a continuously twice differentiable function and $\nabla$ denotes the gradient in $x$-variable. The objective is to obtain the law of the trajectory from its marginals. Let $\Omega = C([0,1], E)$ and ${\mathcal P}(\Omega)$ be the set of probability measures on $\Omega.$  Let $P \in {\mathcal P}(\Omega)$ be the law of $X$ and $\mbW^\sigma$ be the law of $\{\sigma B_t\}_{t \in [0,1]}$ with $B$ being a standard Brownian motion on $E$. They show that (see \cite[Theorem 2.1]{lzks21}) the law  of $X$ can be characterised from its marginals via the following entropy minimisation problem, 
 \begin{equation} \label{eq:enop}
\min \{ H(R \vert \mbW^\sigma) : R \in {\mathcal P}(\Omega), R_t  = P_t \mbox{ for all } t \in  [0,1]\},\end{equation}
where for any $R \in {\mathcal P}(\Omega)$, $R_t$ is the marginal at time $t$.  

The above can be considered for a general diffusion with the generator given by \eqref{eq:cdiff} as discussed in Example \ref{eg:diffusion}. If the associated martingale problem is well-posed, then the one dimensional marginals characterise the law (see \cite[Theorem 4.4.2]{EK}). In addition if hypothesis \eqref{eq:h} and \eqref{eq:dagger} of Theorem \ref{thm:cpmsfm} are satisfied (for e.g.\ when the drift is a bounded continuous function), then  Corollary \ref{cor1} will imply that \eqref{eq:enop} will yield that the minimiser, i.e. $R^\star$ is a Markov process. This has also been observed in \cite[Theorem 4.5]{bl21} assuming uniqueness of solution. 

Further in \cite[Theorem 4.1]{lzks21} it is shown that the $P$ is the unique minimizer using the fact that the Radon-Nikodym derivative of $P$ w.r.t $\mbW^\sigma$ depends only on the marginals. Such an argument will follow in general as long as the quadratic variation process  depends only on the marginals of  $\log \Lambda_t$. This will imply that the minimizer obtained via Corollary \ref{cor1}  yield the unique minimizer as the true law of the process.
\end{example}

We conclude by considering an example of  martingale problems  associated with branching Markov processes.

\begin{example} \label{eg:brandiff}
Let $E_0$ be a Polish space and $\mbu$ be a compact metric space.  Each particle shall   move in $E_0$ according to a  Feller process with generator ${\mathcal B}$, as in Example \ref{eg:diffusion} and \ref{eg:jump}. Each particle branches or dies with a location dependent intensity $\alpha(x,u)$ for $x \in E$ and $u \in U$. We shall assume that $\alpha(\cdot,\cdot)$ is a continuous function on $E_0\times U$. Upon its death it gives rise to children with location dependent offspring distribution whose probability generating function is
$$ \phi(z,x,u) = \sum_{l \geq 0} p_l(x,u) z^l,$$
where $\phi:[0,1]\times \mbR^d \times \mbR^d \mapsto \mbR$ such that $\phi \in C_0([0,1]\times \mbR^d \times \mbR^d)$. 
We assume that the distribution has finite mean , i.e. $\sum_{l \geq 0} l p_l(\cdot,\cdot) < \infty$, and $\sum_{l \geq 0} l p_l(\cdot,\cdot) \in C_0(\mbR^d\times\mbR^d).  $  Let $M_F(E_0)$ denote the space of finite measures on $E_0$ endowed with the topology of weak convergence. For any bounded continuous $f$ on $E_0$, let $$ \langle f , \mu \rangle := \int f d\mu.$$
 Let $\A$  be a linear operator with $\D(\A) = \{ g \in C^2_0(\mbR^d): \parallel g \parallel_\infty <1 \}$  given by
$$ \A g(\mu) =  \exp(\langle \log g, \mu \rangle) \left \langle \frac{{\mathcal B}g + \alpha(\phi(g) -g)}{g}, \mu \right \rangle,$$ 
where $\mu \in M_F(E_0).$ It is easy to see that (A1) and (A2) are satisfied. We can choose  $\B$ as in  Example \ref{eg:diffusion} or \ref{eg:jump}, then by \cite[Theorem 9.4.2]{EK}, (A3) is satisfied. By \cite[Theorem 4.1]{kurtz98} or \cite[Theorem 2.4]{BB}, (A4) holds. Finally, from \cite[Theorem 5.19]{EK},  (A5) holds  when $U$ is as in \eqref{eq:msol} and all parameters are given by bounded continuous functions.  This provides a generic setting where relaxed controlled martingale problems with branching diffusions can have Markov mimics.

 We now turn to an application of Theorem \ref{thm:cpmsfm}. In \cite{bl21},  an  entropy minimization problem with respect to branching Brownian motion is shown to be equivalent to regularized unbalanced optimal transport.
The branching Brownian motion starts with an initial distribution $R_0$ and each particle moves according to a Brownian motion  with diffusion constant $\nu$ in $E = \chi,$ which is a compact smooth Riemannianian manifold without boundary. The branching mechanism is given by ${\bf q} = \{q_k\}_{k \geq 1},$ where $q_k$ is rate at which the particle branches into $k$ particles. We will denote the system of branching Brownian motions by $R\equiv$ BBM($R_0,\nu, {\bf q}$). 

Using stochastic calculus for general semimartingales with jumps, they  show that under exponential moment assumptions on $R_0$ and ${\bf q}$, one can construct modified branching Brownian motions that are  absolutely continuous with BBM($R_0,\nu, {\bf q}$).
In the modified branching Brownian motion, particles move according to a stochastic differential equation with an additional drift $\tilde{v}$  along with time dependent branching rates $\tilde{{\bf q}} = \{\tilde{q}_k(t)\}_{k \geq 1, t \geq 0}$ (see \cite[Theorem 4.23]{bl21} for assumptions on $\tilde{v}$ and  $\tilde{{\bf q}}(t)$).

If one models the trajectory of cell development considered in \cite{lzks21} via a suitable branching diffusion, then an  optimisation problem  with marginal constraints as in \eqref{eq:enop} with ${\mathbb W}^\sigma$ being replaced by $R$ can be considered.  In \cite[Theorem 3.2.1]{ng23}, the optimisation problem 
\begin{equation} \label{eq:enopbr}
\min \{ H(Q \vert R) : Q \in {\mathcal C}\},
\end{equation}
is considered, where
$$
{\mathcal C} = \left \{ Q \in \P(\Omega) : 
                \begin{array}{c}
                Q_t = P_t, \forall t \in [0,1] \mbox{ and } Q \mbox{ is any modified } \\  
                \mbox{branching Brownian motion with branching mechanism }  \tilde{{\bf q}} 
                    \end{array}
   \right \}.
$$
It is shown that the $P$ is the unique minimizer of  \eqref{eq:enopbr} using the fact that the Radon-Nikodym derivative of $P$ w.r.t.\ $R$ depends only on the marginals.


Lastly, if the associated martingale problem for the branching diffusion  is well-posed, then the one dimensional marginals characterise the law (see \cite[Theorem 4.2]{EK}). Thus if the Radon-Nikodym derivative between the branching diffusion and the base branching Brownian motion satisfy hypothesis \eqref{eq:h} and \eqref{eq:dagger} of Theorem \ref{thm:cpmsfm}, (for, e.g.\  when $\tilde{v}$ and  $\tilde{q}(t)$ are bounded continuous), then  Corollary \ref{cor1} will imply that the equivalent problem with \eqref{eq:enopbr} will yield that the minimizer $R^\star$ is a Markov process. One would need additional assumptions as in \cite[Theorem 3.2.1.]{ng23} to show that the unique minimiser is the true law of the process.

\end{example}

%% file: arxivSubmission.bbl
\begin{thebibliography}{10}

\bibitem{BS23}
Sumith~Reddy Anugu and Vivek~S. Borkar.
\newblock A selection procedure for extracting the unique {F}eller weak
  solution of degenerate diffusions.
\newblock {\em Appl. Math. Optim.}, 87(3):46, 2023.

\bibitem{ABG}
Ari Arapostathis, Vivek~S. Borkar, and Mrinal~K. Ghosh.
\newblock {\em Ergodic control of diffusion processes}, volume 143 of {\em
  Encyclopedia of Mathematics and its Applications}.
\newblock Cambridge University Press, Cambridge, 2012.

\bibitem{bF20}
Julio Backhoff-Veraguas and Joaqu{\'\i}n Fontbona.
\newblock Generalized entropy minimization under full marginal constraints.
\newblock {\em arXiv preprint arXiv:2004.10679}, 2020.

\bibitem{bl21}
Aymeric Baradat and Hugo Lavenant.
\newblock Regularized unbalanced optimal transport as entropy minimization with
  respect to branching brownian motion.
\newblock {\em arXiv preprint arXiv:2111.01666}, 2021.

\bibitem{ac20}
Aymeric Baradat and Christian L{\'e}onard.
\newblock Minimizing relative entropy of path measures under marginal
  constraints.
\newblock {\em arXiv preprint arXiv:2001.10920}, 2020.

\bibitem{benes70}
V.~E. Bene\v{s}.
\newblock Existence of optimal strategies based on specified information, for a
  class of stochastic decision problems.
\newblock {\em SIAM J. Control}, 8:179--188, 1970.

\bibitem{BB}
Abhay~G. Bhatt and Vivek~S. Borkar.
\newblock Occupation measures for controlled {M}arkov processes:
  characterization and optimality.
\newblock {\em Ann. Probab.}, 24(3):1531--1562, 1996.

\bibitem{BS10}
V.~S. Borkar and K.~Suresh Kumar.
\newblock A new {M}arkov selection procedure for degenerate diffusions.
\newblock {\em J. Theoret. Probab.}, 23(3):729--747, 2010.

\bibitem{b86}
Vivek~S. Borkar.
\newblock A remark on the attainable distributions of controlled diffusions.
\newblock {\em Stochastics}, 18(1):17--23, 1986.

\bibitem{b91}
Vivek~S. Borkar.
\newblock On extremal solutions to stochastic control problems.
\newblock {\em Appl. Math. Optim.}, 24(3):317--330, 1991.

\bibitem{b93}
Vivek~S. Borkar.
\newblock On extremal solutions to stochastic control problems. {II}.
\newblock {\em Appl. Math. Optim.}, 28(1):49--56, 1993.

\bibitem{bs13}
Gerard Brunick and Steven Shreve.
\newblock Mimicking an {I}t\^{o} process by a solution of a stochastic
  differential equation.
\newblock {\em Ann. Appl. Probab.}, 23(4):1584--1628, 2013.

\bibitem{cgtp21}
Yongxin Chen, Tryphon~T. Georgiou, and Michele Pavon.
\newblock Stochastic control liaisons: {R}ichard {S}inkhorn meets {G}aspard
  {M}onge on a {S}chr\"{o}dinger bridge.
\newblock {\em SIAM Rev.}, 63(2):249--313, 2021.

\bibitem{cgtp22}
Yongxin Chen, Tryphon~T. Georgiou, and Michele Pavon.
\newblock The most likely evolution of diffusing and vanishing particles:
  {S}chr\"{o}dinger bridges with unbalanced marginals.
\newblock {\em SIAM J. Control Optim.}, 60(4):2016--2039, 2022.

\bibitem{pz14}
Giuseppe Da~Prato and Jerzy Zabczyk.
\newblock {\em Stochastic equations in infinite dimensions}, volume 152 of {\em
  Encyclopedia of Mathematics and its Applications}.
\newblock Cambridge University Press, Cambridge, second edition, 2014.

\bibitem{DMNH}
Claude Dellacherie and Paul-Andr\'{e} Meyer.
\newblock {\em Probabilities and potential}.
\newblock North Holland, 1979.

\bibitem{du94}
Bruno Dupire.
\newblock Pricing with a smile.
\newblock {\em Risk}, 7(1):18--20, 1994.

\bibitem{du96}
Bruno Dupire.
\newblock A unified theory of volatility, derivatives pricing: The {C}lassic
  {C}ollection, 2004.

\bibitem{EK}
Stewart~N. Ethier and Thomas~G. Kurtz.
\newblock {\em Markov processes}.
\newblock Wiley Series in Probability and Mathematical Statistics: Probability
  and Mathematical Statistics. John Wiley \& Sons, Inc., New York, 1986.
\newblock Characterization and convergence.

\bibitem{f14}
Martin Forde.
\newblock On the {M}arkovian projection in the {B}runick-{S}hreve mimicking
  result.
\newblock {\em Statist. Probab. Lett.}, 85:98--105, 2014.

\bibitem{ng23}
Nitya Gadhiwala.
\newblock Branching brownian motion models for cell development trajectories.
\newblock Master's thesis, Department of Mathematics, University of British
  Columbia, 2023.
\newblock http://hdl.handle.net/2429/85714.

\bibitem{g86}
I.~Gy\"{o}ngy.
\newblock Mimicking the one-dimensional marginal distributions of processes
  having an {I}t\^{o} differential.
\newblock {\em Probab. Theory Relat. Fields}, 71(4):501--516, 1986.

\bibitem{Kisiel}
Michal Kisielewicz.
\newblock {\em Stochastic differential inclusions and applications}, volume~80
  of {\em Springer Optimization and Its Applications}.
\newblock Springer, New York, 2013.

\bibitem{Krylov}
N.~V. Krylov.
\newblock The selection of a {M}arkov process from a {M}arkov system of
  processes, and the construction of quasidiffusion processes.
\newblock {\em Izv. Akad. Nauk SSSR Ser. Mat.}, 37:691--708, 1973.

\bibitem{Krylov87}
N.~V. Krylov.
\newblock {\em Nonlinear elliptic and parabolic equations of the second order},
  volume~7 of {\em Mathematics and its Applications (Soviet Series)}.
\newblock D. Reidel Publishing Co., Dordrecht, 1987.
\newblock Translated from the Russian by P. L. Buzytsky [P. L. Buzytski\u{\i}].

\bibitem{Krylov08}
N.~V. Krylov.
\newblock {\em Controlled diffusion processes}, volume~14 of {\em Stochastic
  Modelling and Applied Probability}.
\newblock Springer-Verlag, Berlin, 2009.
\newblock Translated from the 1977 Russian original by A. B. Aries, Reprint of
  the 1980 edition.

\bibitem{kurtz98}
Thomas~G. Kurtz and Richard~H. Stockbridge.
\newblock Existence of {M}arkov controls and characterization of optimal
  {M}arkov controls.
\newblock {\em SIAM J. Control Optim.}, 36(2):609--653, 1998.

\bibitem{lzks21}
Hugo Lavenant, Stephen Zhang, Young-Heon Kim, and Geoffrey Schiebinger.
\newblock Towards a mathematical theory of trajectory inference.
\newblock {\em arXiv preprint arXiv:2102.09204}, 2021.

\bibitem{l13}
Christian L{\'e}onard.
\newblock A survey of the {S}chr\"odinger problem and some of its connections
  with optimal transport.
\newblock {\em Discrete \& Continuous Dynamical Systems-A}, 34:1533--1574,
  2014.

\bibitem{lz21}
Qi~L{\"u} and Xu~Zhang.
\newblock Control theory for stochastic distributed parameter systems, an
  engineering perspective.
\newblock {\em Annual Reviews in Control}, 51:268--330, 2021.

\bibitem{m95}
Toshio Mikami.
\newblock Copula fields and their applications.
\newblock {\em Proc. Japan Acad. Ser. A Math. Sci.}, 71(10):221--224 (1996),
  1995.

\bibitem{m99}
Toshio Mikami.
\newblock Markov marginal problems and their applications to {M}arkov optimal
  control.
\newblock In {\em Stochastic analysis, control, optimization and applications
  (W. McEneaney, G. Yin, Q. Zhang, eds.)}, Systems Control Found. Appl., pages
  457--476. Birkh\"{a}user Boston, Boston, MA, 1999.

\bibitem{m21}
Toshio Mikami.
\newblock {\em Stochastic optimal transportation---stochastic control with
  fixed marginals}.
\newblock Springer Briefs in Mathematics. Springer, Singapore, [2021]
  \copyright 2021.

\bibitem{SV}
Daniel~W. Stroock and S.~R.~Srinivasa Varadhan.
\newblock {\em Multidimensional diffusion processes}.
\newblock Classics in Mathematics. Springer-Verlag, Berlin, 2006.
\newblock Reprint of the 1997 edition.

\bibitem{yos95}
K\'osaku Yosida.
\newblock {\em Functional analysis}.
\newblock Classics in Mathematics. Springer-Verlag, Berlin, 1995.
\newblock Reprint of the sixth (1980) edition.

\end{thebibliography}
